\documentclass[letterpaper, 10 pt, conference]{ieeeconf}
\IEEEoverridecommandlockouts        

\usepackage{cite}
\usepackage{amsmath}
\usepackage{amssymb,amsfonts}
\usepackage{algorithmic}
\usepackage{hyperref}
\usepackage{textcomp}
\usepackage{graphicx,color}
\usepackage{mathrsfs}
\usepackage[vlined,ruled]{algorithm2e}
\usepackage{subfigure}
\usepackage{url}
\usepackage{color}
\usepackage{dsfont}
\usepackage{bbm}
\usepackage{booktabs}
\usepackage{array}
\usepackage[table]{xcolor} 
\usepackage{yfonts}
\usepackage{float}
\usepackage{afterpage}
\usepackage{mathabx}

\newtheorem{theorem}{Theorem}[section]
\newtheorem{lemma}[theorem]{Lemma}

\newcommand{\real}{\mathbb{R}}

\DeclareSymbolFont{bbold}{U}{bbold}{m}{n}
\DeclareSymbolFontAlphabet{\mathbbold}{bbold}

\newcommand\oprocendsymbol{\hbox{$\square$}}
\newcommand\oprocend{\relax\ifmmode\else\unskip\hfill\fi\oprocendsymbol}

\renewcommand{\theenumi}{(\roman{enumi})}

\newcommand*{\QEDA}{\hfill\ensuremath{\blacksquare}}%

\makeatletter
\renewcommand*{\@opargbegintheorem}[3]{\trivlist
      \item[\hskip \labelsep{\bfseries #1\ #2}] \textbf{(#3)}\ \itshape}
\makeatother

\graphicspath{{img/}}

\begin{document}
\title{\bf On Direct vs Indirect Data-Driven Predictive
  Control} \author{Vishaal Krishnan \qquad Fabio Pasqualetti
  \thanks{This material is based upon work supported in part by awards
    AFOSR FA9550-20-1-0140 and FA9550-19-1-0235, and ONR
    N00014-19-1-2264. The authors are with the Department of
    Mechanical Engineering, University of California at Riverside,
    Riverside, CA 92521 USA. E-mail:
    \href{mailto:vishaalk@ucr.edu}{\texttt{vishaalk@ucr.edu}},
    \href{mailto:fabiopas@engr.ucr.edu}{\texttt{fabiopas@engr.ucr.edu.}}
  } }
\maketitle

\thispagestyle{empty}

\begin{abstract}
In this work, we compare the direct and indirect approaches
to data-driven predictive control of stochastic linear
time-invariant systems. The distinction between the two approaches lies in the fact
that the indirect approach involves identifying a lower dimensional 
model from data which is then used in a certainty-equivalent control design, 
while the direct approach avoids this intermediate step altogether.
Working within an optimization-based framework,
we find that the suboptimality gap measuring the control performance 
w.r.t. the optimal model-based control design vanishes with the size of the dataset only with the direct approach. 
The indirect approach has a higher rate of
convergence, but its suboptimality gap does not vanish
as the size of the dataset increases. 
This reveals the existence of two distinct regimes of performance
as the size of the dataset of input-output behaviors is increased.
We show that the indirect approach, by relying on the identification 
of a lower dimensional model, has lower variance and 
outperforms the direct approach for smaller datasets, 
while it incurs an asymptotic bias as a result of the process noise 
and a (possibly) incorrect assumption on the order of the identified model. 
The direct approach, however, does not incur an asymptotic bias, and outperforms 
the indirect approach for larger datasets. 
Ultimately, by revealing the existence of two non-asymptotic regimes 
for the performance of direct and indirect data-driven predictive
control designs, our study suggests that neither approach is invariably 
superior and that the choice of design must, in practice, be informed 
by the available dataset.
\end{abstract}


\section{Introduction}
\label{sec:intro}
The interest in direct data-driven control is motivated
by its promise to make the system identification step 
unnecessary for control design, and more optimistically by its bid to outperform 
the traditional system identification and model-based control design pipeline.
While there has been a surge in literature devoted to developing
techniques for direct data-driven control 
\cite{CDP-PT:19,GB-VK-FP:19,JB-AK-CWS-FA:20,HJV-MKC-MM:20,GB-DSB-FP:20}, 
the question of its success on the above fronts has largely remained unsettled. 
With this broad motivation, we seek to systematically compare 
the performance of the two paradigms in the predictive 
control setting for stochastic linear time-invariant systems.  
 
We consider a discrete-time, stochastic
linear time-invariant (LTI) system of the form: 
\begin{align}
	\begin{aligned}
		x_{t+1} = A x_t + B u_t + w_t, \qquad
		y_t = C x_t,
	\end{aligned}
	\label{eq:LTI_system}
\end{align}
where $A \in \real^{n \times n}$, $B \in \real^{n \times m}$
and $C \in \real^{p \times n}$ are the system matrices, 
$x_t \in \real^n$, $u_t \in \real^m$ and $y_t \in \real^p$
are the system state, control input and output at 
time~$t \in \mathbb{N}$, respectively,
with $w_t \in \real^n$ being the process noise
generated by an i.i.d. process with distribution 
$\mathcal{N}(0, \Omega_{\rm w})$.
Furthermore, we assume that the pair $(A,C)$ is 
observable and that $(A,B)$ is controllable.
The predictive control task is specified via the following optimization problem:
\begin{align} \label{eq:ctrl_task}
\begin{aligned}
 	\min_{u_0, \ldots, u_{T-1}}	~&\mathbb{E} \left[ \sum_{t=0}^{T-1} \left( \alpha_t(u_t) + \beta_t(y_{t+1}) \right) \right] \\
 				&\text{s.t.} \quad \begin{cases} y_t &= C x_t, \\ x_{t+1} &= A x_t + B u_t + w_t, \\ x_0 &= 0,  \end{cases}
 \end{aligned}
\end{align}
where the stage costs $\alpha_t(u) = u^\top Q_t u$ and $\beta_t(y) = (y - y_{\rm ref})^\top R_t (y - y_{\rm ref})$,
with $Q_t > 0$ and $R_t > 0$ for all $t \in \lbrace 0, \ldots, T-1 \rbrace$
and $y_{\rm ref} \in \real^p$.
The data-driven control problem studied in this paper is one of solving 
the control task~\eqref{eq:ctrl_task} using a finite dataset of input-output
behaviors of System~\eqref{eq:LTI_system}.
To this end, we assume that the system matrices~$A$, $B$ 
and $C$, and the process noise covariance $\Omega_{\rm w}$ 
are unknown. Instead, we have access to
the input-output data from $N$ control experiments over
a time horizon~$\lbrace 0, \ldots, T \rbrace$ on 
System~\eqref{eq:LTI_system} with initial condition~$x_0 = 0$:
\begin{align*}
	U = \left[ \begin{matrix} \mathbf{u}^{(1)} & \ldots & \mathbf{u}^{(N)} \end{matrix} \right], \quad
	Y = \left[ \begin{matrix} \mathbf{y}^{(1)} & \ldots & \mathbf{y}^{(N)} \end{matrix} \right],
\end{align*}
where $\mathbf{u}^{(i)} = (u_0^{(i)}, \ldots, u_{T-1}^{(i)})$ 
and $\mathbf{y}^{(i)} = (y_1^{(i)}, \ldots, y_T^{(i)})$
for every $i \in \lbrace 1, \ldots, N \rbrace$.
We assume that the inputs $u_t^{(i)}$ for the control 
experiments are generated by an i.i.d. process with 
distribution $\mathcal{N}(0, \Sigma_{\rm u})$.

Data-driven predictive control design essentially
involves the mapping of the dataset of control experiments 
and the parameters of the control task
onto a finite control sequence.
Optimality of the control design depends on extracting
system and task-relevant information from the 
dataset~\cite{AA-SS:18,VP-AM:19,NT-FP-WB:00}. The presence of noise in the dataset introduces system 
and task-irrelevant information, and since both system 
and task are unknown, extracting relevant information 
from the dataset is key. If the information
extraction step is not handled effectively,
it results in suboptimality of the control design.
Against this backdrop, the distinction between 
the direct and indirect approaches to data-driven control
lie in the way the two approaches extract system and task-relevant
information from the dataset.
Yet, a complete understanding of the comparative
advantages of the two approaches is lacking.
Motivated by this need, we undertake in this paper
a comparative study between the two approaches.
We choose for comparison with the direct data-driven 
control design a certainty equivalent control design 
that utilizes a model identified from the dataset 
by ordinary least squares. We later discuss the
implications of this choice for our results and
conclusions, and emphasize that the 
qualitative insights we obtain on the 
comparative performances are much more general
and not constrained by this choice. 

\textbf{Contributions.}
The primary contribution of this paper is a
comparative study of the non-asymptotic 
performance of direct and indirect predictive control
for stochastic LTI systems. Our investigation reveals the existence of
two distinct non-asymptotic regimes for performance
as the dataset size is increased,
where one approach outperforms the
other in each regime. More specifically,
we find that the indirect approach, by relying
on the identification of a low dimensional
kernel representation, converges faster to its
asymptotic performance (measured by its suboptimality gap) 
and outperforms the direct approach for smaller datasets.
Conversely, the direct approach outperforms the indirect
approach for larger datasets and achieves better
asymptotic performance with a vanishing suboptimality gap.
Our technical contributions are as follows:
We first formulate a unifying optimization-based framework that
permits a comparative analysis of the direct and indirect 
approaches. We then present analytical results characterizing the
asymptotic performance and sample complexity 
bounds for the two approaches, shedding light on their
performances as the dataset size is increased, and their
dependence on the system and noise parameters.
Finally, we present results from numerical experiments
validating our analysis and demonstrating the existence
of the aforementioned non-asymptotic regimes in practice.

\textbf{Related work.}
The control design problem considered in this paper, 
that of data-driven predictive control, has garnered considerable
attention in recent years. Several techniques for design
that combine learning with model predictive control
have been proposed~\cite{LH-KW-MM-MZ:20,JB-JK-MAM-FA:20}. In~\cite{JC-JL-FD:18}, 
the authors introduce the Data-Enabled Predictive Control 
(DeePC) method, which has then been applied to
various settings~\cite{JC-JL-FD:20}.
In a recent work~\cite{FD-JC-IM:21},
the authors present experimental results comparing direct and
indirect data-driven control in the presence of noise and
nonlinearities in the underlying system. 
Furthermore, the authors present a framework to bridge
the direct and indirect approaches based on regularization,
which potentially allows for an efficient transition between the
two approaches. Also, in~\cite{FF-SL:21} the authors compare the DeePC
method with the Subspace Predictive Control (SPC) method,
showing that the two methods are equivalent in the
deterministic case, reasoning that the DeePC method
implicitly estimates the same predictive model as SPC.
Further, they investigate the comparative performance
of the two methods in experiments. Also of particular relevance are works that
have investigated the connection between persistency of excitation
and data-driven control and system identification~\cite{HJVW-JE-HLT-MKC:20,JCW-PR-IM-BLMDM:05},
as are works that have investigated the sample complexity 
of LTI system identification~\cite{MS-HM-ST-MIJ-BR:18,YJ-AP:19} 
and data-driven LQR design~\cite{SD-HM-NM-BR-ST:18}.
However, an analysis characterizing the 
comparative performance of the direct and indirect approaches
and an understanding of the non-asymptotic case are still lacking.

The distinction between direct and indirect
data-driven control parallels the distinction between 
discriminative and generative modeling approaches 
in machine learning classification, where a dataset is either used to 
train a classifier to simply learn (in a supervised
setting) to classify input examples as in the discriminative case,
or to learn a generative model for the classes as in the latter case.
In this context, \cite{AN-MIJ:02} 
is an early work comparing the performance
of discriminative and generative classifiers as a function of dataset size.
Several works~\cite{GB-BT:04,JL-CB-TM:06} have proposed interpolating schemes that 
result in hybrid discriminative/generative models to exploit the
advantages of both classes of models, with reported success.

\section{Data-driven control design}
In this section, we address the data-driven control design problem.
To this end, we first take the view that the data-driven predictive 
control design problem is one of solving the control task~\eqref{eq:ctrl_task}
with a data-driven model of System~\eqref{eq:LTI_system}. We then 
note that both the direct and indirect approaches can be studied
within this framework by showing that even the direct approach
to data-driven control design relies on an implicit model of the input-output 
behavior of the underlying system, even if such a model is not explicitly identified
(as recently noted also in~\cite{FF-SL:21}). 
We then obtain a characterization of the suboptimality gap for
the data-driven formulation of the control task~\eqref{eq:ctrl_task}, 
showing that it is controlled both above and below by the 
discrepancy between the model used in design and the 
true model. This allows for a comparison of the direct and indirect approaches 
to data-driven control, via the discrepancy between the implict and true models in the 
case of the direct approach, and between the identified and true models for 
the indirect approach, respectively.
 
We begin by establishing a key property of the underlying 
control task~\eqref{eq:ctrl_task} in the model-based setting,
which will be useful in setting up the data-driven formulation.
The output $\mathbf{y} = (y_1, \ldots, y_T)$ of System~\eqref{eq:LTI_system}  
(note that we have $y_0 = 0$ since $x_0 = 0$)
over the horizon $\lbrace 1, \ldots, T \rbrace$
generated by control input $\mathbf{u}
= (u_0, \ldots, u_{T-1})$ is given by:
\begin{align} \label{eq:IO_response}
	\mathbf{y} = \mathcal{G} \mathbf{u} + \mathcal{G}' \mathbf{w},
\end{align}
where $\mathbf{w} = (w_0, \ldots, w_{T-1})$ is the process noise, and:
\begin{align*}
	\mathcal{G} =  \left[ \begin{matrix} 	CB & 0 & \ldots & 0 \\
																CAB & CB & \ldots & 0 \\
																	\vdots & \vdots & \ddots & \vdots \\
																		C A^{T-1} B & C A^{T-2} B & \ldots & CB \end{matrix} \right].
\end{align*}
$\mathcal{G}'$ is obtained by replacing $B$ with $I_n$ in the expression for $\mathcal{G}$ above.
Now, with $\mathbf{v} = \mathcal{G}' \mathbf{w}$ in~\eqref{eq:IO_response}
and  $F(\mathbf{u}, \mathbf{y}) = \sum_{t=0}^{T-1} \left( \alpha_t(u_t) + \beta_t(y_{t+1}) \right)$,
we can express the control task~\eqref{eq:ctrl_task} as follows:
\begin{align} \label{eq:behavioral_formulation_ctrl_task}
	\min_{\mathbf{u} \in \real^{mT}}~ \mathbb{E}_{\mathbf{v}} \left \lbrace F(\mathbf{u}, \mathbf{y}) \; | \; \mathbf{y} = \mathcal{G} \mathbf{u} + \mathbf{v} \right \rbrace.
\end{align}
The following lemma establishes the certainty equivalence 
property for model-based predictive control of stochastic 
LTI systems with quadratic cost:
\begin{lemma}[\bf \emph{Certainty equivalent model-based predictive control}] \label{lemma:certainty_equiv_model-based}
For $\mathbf{u}^* = \arg \min_{\mathbf{u} \in \real^{mT}}~ F(\mathbf{u}, \mathcal{G} \mathbf{u} )$,
we have $\mathbb{E}_{\mathbf{v}} \left[ F(\mathbf{u}^*, \mathcal{G} \mathbf{u}^* + \mathbf{v}) \right] 
= \min_{\mathbf{u} \in \real^{mT}} \mathbb{E}_{\mathbf{v}} \left[ F(\mathbf{u}, \mathcal{G} \mathbf{u} + \mathbf{v}) \right]$. \oprocend
\end{lemma}
\medskip
\noindent
We refer the reader to Appendix~\ref{app:proof_certainty_equiv} for the proof.
Lemma~\ref{lemma:certainty_equiv_model-based} establishes that
given the input-output model~$\mathcal{G}$, certainty equivalence holds
w.r.t. process noise. Therefore, the minimizer to
the stochastic optimization problem~\eqref{eq:behavioral_formulation_ctrl_task}
can be equivalently obtained as the solution to the following
deterministic optimization problem\footnote{We use the notation $F(\mathbf{u}, \mathbf{y})$ and $F(\mathcal{P}\mathbf{u})$
interchangeably. We see that $\mathcal{P}\mathbf{u} = (\mathbf{u}, \mathcal{G}\mathbf{u}) \in \real^{(m+p)T}$
is indeed of the same dimension as $(\mathbf{u}, \mathbf{y}) \in \real^{(m+p)T}$.}: 
\begin{align} \label{eq:ctrl_design_problem}
	\min_{\mathbf{u} \in \real^{mT}}~ F(\mathcal{P} \mathbf{u}), \qquad \mathcal{P} = \left[ \begin{matrix} I \\ \mathcal{G} \end{matrix} \right].
\end{align}
This suggests the use of the above formulation~\eqref{eq:ctrl_design_problem}
as the control design procedure for the task~\eqref{eq:ctrl_task}.
We now compute the minimizer~$\mathbf{u}^*$ in~\eqref{eq:ctrl_design_problem}. 
With $\mathbf{y}_{\rm ref} = I_T \otimes y_{\rm ref}$, we get:
\begin{align*}
	F(\mathbf{u}, \mathbf{y}) = \mathbf{u}^\top \mathcal{Q} \mathbf{u} 
												+ \left( \mathbf{y} - \mathbf{y}_{\rm ref} \right)^\top \mathcal{R} \left( \mathbf{y} - \mathbf{y}_{\rm ref} \right),
\end{align*}
where $\mathcal{Q} = \mathrm{diag} \left( Q_0, \ldots, Q_{T-1} \right)$
and $\mathcal{R} = \mathrm{diag} \left( R_0, \ldots, R_{T-1} \right)$.
We further obtain:
\begin{align*}
	\nabla F(\mathbf{u}, \mathbf{y}) = 2 \left[ \begin{matrix} \mathcal{Q} & 0 \\ 0 & \mathcal{R} \end{matrix} \right] 
																		\left( \begin{matrix} \mathbf{u} \\ \mathbf{y} \end{matrix} \right)  
																	- 2 \left( \begin{matrix} 0 \\ \mathcal{R} \mathbf{y}_{\rm ref} \end{matrix} \right).
\end{align*}
The minimizer $\mathbf{u}^*$ in~\eqref{eq:ctrl_design_problem}
satisfies $\mathcal{P}^\top \nabla F (\mathcal{P} \mathbf{u}^*) = 0$
(first-order optimality condition).
Substituting from the above, we get:
\begin{align*}
	\mathcal{P}^\top \left[ \begin{matrix} \mathcal{Q} & 0 \\ 0 & \mathcal{R} \end{matrix} \right] \mathcal{P} \mathbf{u}^*
			= \mathcal{P}^\top \left( \begin{matrix} 0 \\ \mathcal{R} \mathbf{y}_{\rm ref} \end{matrix} \right),
\end{align*}
Simplifying the above, we get:
\begin{align} \label{eq:optimal_ctrl_input_model-based}
	\mathbf{u}^* = \left( \mathcal{Q} + \mathcal{G}^\top \mathcal{R} \mathcal{G} \right)^{-1} \mathcal{G}^\top \mathcal{R} \mathbf{y}_{\rm ref}. 
\end{align}
We note however, that the certainty equivalence property 
was established in Lemma~\ref{lemma:certainty_equiv_model-based} 
under the assumption of availability of the true 
input-output behavior model~$\mathcal{P}$.
In the data-driven control setting, we do not have direct access to
$\mathcal{P}$ for control design by~\eqref{eq:ctrl_design_problem},
and either (i) an estimate of $\mathcal{P}$ is obtained from 
noisy input-output behavior data $U, Y$ (indirect data-driven control),
or (ii) the data matrix $\left[ \begin{matrix} U^\top & Y^\top \end{matrix} \right]^\top$ 
is itself used in place of the behavior model (direct data-driven control).
In other words, data-driven control design involves the use
of an estimate $\widehat{\mathcal{P}}$ of the true behavior 
model $\mathcal{P}$. 
In indirect data-driven control design,
such an estimate is explicitly obtained from data, whereas 
direct data-driven control design involves the use of an implicit
estimate, as will be seen in the ensuing section. The
data-driven control design problem is formulated
by replacing $\widehat{\mathcal{P}}$ for $\mathcal{P}$ 
in~\eqref{eq:ctrl_design_problem}:
\begin{align} \label{eq:data-driven_ctrl_design}
	 \min_{\mathbf{u} \in \real^{mT}} F \left(\widehat{\mathcal{P}} \mathbf{u} \right), \qquad \widehat{\mathcal{P}} = \left[ \begin{matrix} I \\ \widehat{\mathcal{G}} \end{matrix} \right].
\end{align}
Also, let $\widehat{\mathbf{u}} = \arg \min_{\mathbf{u} \in \real^{mT}} F(\widehat{\mathcal{P}} \mathbf{u})$.
Following similar steps as in the computation of $\mathbf{u}^*$
in \eqref{eq:optimal_ctrl_input_model-based}, we get:
\begin{align} \label{eq:optimal_ctrl_input_estimate}
	\widehat{\mathbf{u}} 
		= \left( \mathcal{Q} + \widehat{\mathcal{G}}^\top \mathcal{R} \widehat{\mathcal{G}} \right)^{-1} \widehat{\mathcal{G}}^\top \mathcal{R} \mathbf{y}_{\rm ref}. 
\end{align}
We note that the control input $\widehat{\mathbf{u}}$
is not guaranteed to be optimal for the control task~\eqref{eq:ctrl_task},
owing to the mismatch between the estimate $\widehat{\mathcal{G}}$
and the true model~$\mathcal{G}$. The performance of $\widehat{\mathbf{u}}$
is measured by its suboptimality gap given~by:
\begin{align} \label{eq:suboptimality_gap}
	{\rm Gap} (\widehat{\mathbf{u}}) = F \left(\mathcal{P} \widehat{\mathbf{u}} \right) - F(\mathcal{P} \mathbf{u}^*).
\end{align}
This suboptimality gap serves as a metric for comparing the direct and indirect data-driven 
control design methodologies.
We note that ${\rm Gap} (\widehat{\mathbf{u}}) \geq 0$ for any $\widehat{\mathbf{u}}$, and 
from the $\mu$-strong convexity of $F$\footnote{where $\mu =  \lambda_{\rm min} \left( \mathcal{Q} + \mathcal{G}^\top \mathcal{R} \mathcal{G} \right)$
and $\nu =  \lambda_{\rm max} \left( \mathcal{Q} + \mathcal{G}^\top \mathcal{R} \mathcal{G} \right)$. \label{footnote:func_F}}, 
it follows that ${\rm Gap} (\widehat{\mathbf{u}}) \geq \mu \left \| \left( \widehat{\mathbf{u}} - \mathbf{u}^* \right) \right \|^2 / 2$.
From the $\nu$-Lipschitz continuity of the gradient of~$F^{\ref{footnote:func_F}}$, 
it follows that ${\rm Gap} (\widehat{\mathbf{u}}) \leq \nu \left \| \left( \widehat{\mathbf{u}} - \mathbf{u}^* \right) \right \|^2 / 2$.
Combining the above, we get:
\begin{align*}
  	\frac{\mu}{2} \left \| \left( \widehat{\mathbf{u}} - \mathbf{u}^* \right) \right \|^2
  							\leq {\rm Gap} (\widehat{\mathbf{u}}) 
  							\leq \frac{\nu}{2} \left \| \left( \widehat{\mathbf{u}} - \mathbf{u}^* \right) \right \|^2.
\end{align*}
It follows from~\eqref{eq:optimal_ctrl_input_model-based} 
and~\eqref{eq:optimal_ctrl_input_estimate} that
the error $\| \widehat{\mathbf{u}} - \mathbf{u}^* \|$ 
arises from the mismatch~$\Delta$ between 
$\widehat{\mathcal{G}}$ and $\mathcal{G}$,
and we see from the above that it also controls the suboptimality
gap ${\rm Gap}(\widehat{\mathbf{u}})$.
This allows us to investigate the control performance 
measured by the suboptimality gap ${\rm Gap}(\widehat{\mathbf{u}})$
via the mismatch~$\Delta$ between~$\widehat{\mathcal{G}}$ and $\mathcal{G}$\footnote{It is 
further possible to obtain a bound on $\| \widehat{\mathbf{u}} - \mathbf{u}^* \|$
as a function of $\Delta$ and the system and task parameters. We do not,
however, pursue a characterization of such a bound in the context of this paper.}. 
\subsection{Direct data-driven control}
In direct data-driven control design, we would like to use the 
input-output behavior data matrix $\left[ \begin{matrix} U^\top & Y^\top \end{matrix} \right]^\top$
directly in place of the model $\widehat{\mathcal{P}}$ in~\eqref{eq:data-driven_ctrl_design}.
The key idea here is to avoid identifying a model of input-output behaviors
and to directly search for the optimal behavior for the task~\eqref{eq:ctrl_task} 
within the span of observed behaviors contained in the 
data matrix $\left[ \begin{matrix} U^\top & Y^\top \end{matrix} \right]^\top$.
However, before obtaining the direct data-driven design formulation, 
we first note that the input-output behavior data satisfies:
\begin{align*}
	\left[ \begin{matrix} U \\ Y \end{matrix} \right] = \mathcal{P} U + \left[ \begin{matrix} 0 \\ V \end{matrix}	 \right],
\end{align*}
where $V = \left[ \begin{matrix} \mathbf{v}^{(1)} & \ldots & \mathbf{v}^{(N)}	\end{matrix} \right]$
is the matrix of noise realizations in the control experiments. 
We note that while the above relation reveals the underlying structure in the
available input-output behavior data, we do not have access to $\mathcal{P}$ and~$V$.
The above relation reveals that as the number of experiments~$N$
increases, it may be possible to construct behaviors
$\left[ \begin{matrix} U^\top & Y^\top \end{matrix} \right]^\top \mathbf{z}$,
with $\mathbf{z} \in \real^N$, such that $\mathbf{z} \in \mathrm{Ker}(U)$,
but which nevertheless incur a low cost in design, i.e., 
$F(\left[ \begin{matrix} U^\top & Y^\top \end{matrix} \right]^\top \mathbf{z})$
attains a low value. The corresponding control input would indeed be 
$\mathbf{u} = U \mathbf{z} = 0$, suggesting that the low cost 
can be attained without taking any control action.
However, this is entirely misleading, as such behaviors are
essentially constructed from the process noise components realized in
the control experiments and contained in~$V$, and the corresponding
input $\mathbf{u} = 0$ may not actually incur a low suboptimality gap,
as measured by~${\rm Gap}(\mathbf{u})$.
We therefore restrict our search within $\mathrm{Ker}^\perp (U)$,
the orthogonal complement of~$\mathrm{Ker}(U)$, to
obtain the direct data-driven design formulation:
\begin{align} \label{eq:direct_data-driven_design}
	\min_{\mathbf{z} \in \real^N}~ F \left( \left[ \begin{matrix} U \\ Y\end{matrix} \right]  U^\dagger U \mathbf{z} \right).
\end{align}
Now, for any $\mathbf{z}' = U^\dagger U \mathbf{z} \in \mathrm{Ker}^\perp(U)$, we have:
\begin{align*}
	&\left[ \begin{matrix} U \\ Y \end{matrix} \right] \mathbf{z}' 
			= \left[ \begin{matrix} U \\ Y \end{matrix} \right] U^\dagger U \mathbf{z} 
			= \left( \mathcal{P} U + \left[ \begin{matrix} 0 \\ V \end{matrix}	 \right] \right) U^\dagger U \mathbf{z} \\
			&= \left( \mathcal{P} U U^\dagger U  + \left[ \begin{matrix} 0 \\ V U^\dagger U \end{matrix}	 \right] \right)  \mathbf{z}
			= \left( \mathcal{P} U + \left[ \begin{matrix} 0 \\ V U^\dagger U \end{matrix}	 \right] \right)  \mathbf{z} \\
			&= \left( \mathcal{P} + \left[ \begin{matrix} 0 \\ V U^\dagger \end{matrix}	 \right] \right)  U \mathbf{z} 
			 = \left( \mathcal{P} + \left[ \begin{matrix} 0 \\ V U^\dagger \end{matrix}	 \right] \right) \mathbf{u},
\end{align*}
where $\mathbf{u} = U \mathbf{z} \in \mathrm{Col}(U)$.
With $\widehat{P}_{\rm direct} = \mathcal{P} + \left[ \begin{matrix} 0 \\ V U^\dagger \end{matrix}	 \right]$,
we can rewrite~\eqref{eq:direct_data-driven_design} as:
\begin{align} \label{eq:direct_data-driven_implicit}
	\min_{\mathbf{u} \in \mathrm{Col}(U)} ~ F \left( \widehat{P}_{\rm direct} \mathbf{u} \right).
\end{align}
The above is the sense in which the direct data-driven control
design formulation employs implicitly an estimate $\widehat{P}_{\rm direct}$ of the
true behavior model $\mathcal{P}$ as stated earlier, and can be
connected to the general data-driven control design formulation~\eqref{eq:data-driven_ctrl_design}.
We now obtain the minimizer $\widehat{\mathbf{u}}_{\rm direct}$
in~\eqref{eq:direct_data-driven_design}~as:
\begin{align*}
	\widehat{\mathbf{u}}_{\rm direct} = U \left( U^\top \mathcal{Q} U + Y_{|U}^\top \mathcal{R} Y_{|U} \right)^{\dagger} Y_{|U}^\top \mathcal{R} \mathbf{y}_{\rm ref},
\end{align*}
where $Y_{|U} = Y U^\dagger U$.
Furthermore, we note that the mismatch between $\widehat{\mathcal{G}}_{\rm direct}$
and $\mathcal{G}$ is given by: 
\begin{align} \label{eq:implicit_model_mismatch}
	\Delta_{\rm direct} = \widehat{\mathcal{G}}_{\rm direct} - \mathcal{G} = V U^\dagger.
\end{align}
We characterize the dependence of the implicit model error 
on the number of control experiments~$N$ and the time horizon~$T$
through the following theorem, for the Single-Input-Single-Output case ($p=m=1$) 
for the sake of simplicity. We note that the result 
can be readily extended to the Multiple-Input-Multiple-Output case.
\begin{theorem}[\bf \emph{Implicit model}] \label{thm:direct_data-driven}
Let $p = m =1$ and $N, T \in \mathbb{N}$ be such that
the empirical covariance matrix
$\Sigma_{\rm uu}^N = U U^\top / N$ is invertible.
The implicit model error $\Delta_{\rm direct}$, 
given by~\eqref{eq:implicit_model_mismatch}, satisfies:
\begin{align*}
	&\textit{(i)}	\quad \mathbb{E} \left[ \Delta_{\rm direct} \right] = 0, \\	
	&\textit{(ii)} \quad \mathbb{P} \left\lbrace \left\| \Delta_{\rm direct} \right\|_F \geq \epsilon \right\rbrace 
			\leq \frac{T^2}{N \epsilon^2} \cdot \frac{ \sigma_{\rm w} \sigma_{\rm u}}{\sigma_{\min}^2 \left( \Sigma_{\rm uu}^N \right)},
\end{align*}
where $\sigma_{\min} \left( \Sigma_{\rm uu}^N \right)$ is the smallest
singular value of $\Sigma_{\rm uu}^N$, $\sigma_{\rm u}$ the
variance of input $u$ and $\sigma_{\rm w} 
= \sum_{t=0}^{T-1}  CA^t \Omega_{\rm w} A^{t \top} C^\top$. \oprocend
\end{theorem}
\medskip
\noindent
We refer the reader to Appendix~\ref{app:proof_direct} for the proof.
Theorem~\ref{thm:direct_data-driven} sheds light on both the 
asymptotic and non-asymptotic performance of the direct approach
as the dataset size~$N$ increases.
It establishes that the direct approach does not incur an 
asymptotic bias, which converges to zero ${\rm w.p.}~1$ at 
the rate $\mathcal{O}\left( \frac{T^2}{N} \right)$. The~$T^2$
dependence has implications for the scalability of 
performance of the direct approach which, as we will also see
from numerical experiments, deteriorates drastically with the control horizon
length~$T$.
\subsection{Indirect data-driven control}
For the indirect approach, as stated earlier, we choose as the
candidate a certainty-equivalent control design that utilizes a model
identified from the dataset by ordinary least squares. To develop the
indirect data-driven formulation, we first note
that~\eqref{eq:IO_response} can be rewritten as
$\left[ \begin{matrix} - \mathcal{G} & I \end{matrix} \right]
\left( \begin{matrix} \mathbf{u} \\ \mathbf{y} -
    \mathbf{v} \end{matrix} \right) = 0$.  Let
$\mathcal{M}_L \in \real^{pT \times pT}$ be a lower block-triangular
block-Toeplitz matrix with first column block
$\left[ \begin{matrix} M_1^\top & \ldots & M_{L}^\top & 0 & \ldots &
    0 \end{matrix} \right]^\top$ and row block
$\left[ \begin{matrix} M_1 & 0 & \ldots & 0 \end{matrix} \right]$,
with $M_k \in \real^{p \times p}$ for all
$p \in \lbrace 1, \ldots, L \rbrace$.  For $L \geq n+1$, we note that
there exists $M_k \in \real^{p \times p}$,
$p \in \lbrace 1, \ldots, L \rbrace$, with $M_1 = I_p$ such that the
matrix
$\mathcal{N}_L = - \mathcal{M}_L \mathcal{G} \in \real^{pT \times mT}$
is also lower block-triangular block-Toeplitz, with first column block
$\left[ \begin{matrix} N_1^\top & \ldots & N_L^\top & 0 & \ldots &
    0 \end{matrix} \right]^\top$ and row block
$\left[ \begin{matrix} N_1 & 0 & \ldots & 0 \end{matrix} \right]$,
with $N_k \in \real^{p \times m}$ for all
$p \in \lbrace 1, \ldots, L \rbrace$ and $N_L = 0$.\footnote{This
  follows from the fact that the observability matrix generated by the
  pair $(A,C)$ attains full column rank over a horizon of length~$n$.}
We then get:
\begin{align} \label{eq:IO_delay_rep_vec_form}
\small
\begin{aligned}
	 \left[ \begin{matrix} \mathcal{N}_L & \mathcal{M}_L \end{matrix} \right]
		\left( \begin{matrix} \mathbf{u} \\ \mathbf{y} - \mathbf{v} \end{matrix} \right)
		= \mathcal{M}_L \left[ \begin{matrix} - \mathcal{G} & I \end{matrix} \right]
		\left( \begin{matrix} \mathbf{u} \\ \mathbf{y} - \mathbf{v} \end{matrix} \right)
		= 0.
\end{aligned}
\end{align}
It can be readily seen that this corresponds to the delay operator representation~\cite{CDP-PT:19},
with $L \geq n+1$:
\begin{align} \label{eq:IO_delay_rep}
\begin{aligned}
	\sum_{\tau =1}^L M_{\tau} \left( y_{t-\tau+1} - v_{t-\tau+1} \right) + \sum_{\tau =1}^L N_{\tau} u_{t-\tau} = 0.
\end{aligned}
\end{align}
Now since a pair $\left( \mathcal{M}_L, \mathcal{N}_L \right)$ 
satisfying $\mathcal{N}_L = - \mathcal{M}_L \mathcal{G}$ 
exists for any $L \geq n+1$, we have
$\left[ \begin{matrix} \mathcal{N}_L & \mathcal{M}_L \end{matrix} \right] \mathcal{P}
= \left[ \begin{matrix} \mathcal{N}_L & \mathcal{M}_L \end{matrix} \right] 
\left[ \begin{matrix} I \\ \mathcal{G} \end{matrix} \right] = 0$.
Furthermore, since $M_1 = I_p$ we note that $\mathcal{M}_L$ is invertible for any $L$.
Therefore, $\mathcal{G} = - \mathcal{M}_L^{-1} \mathcal{N}_L$,
and any such pair $(\mathcal{M}_L, \mathcal{N}_L)$ (for $L \geq n+1$)
satisfying the above allows for an equivalent representation of the 
input-output behavior model~$\mathcal{P}$. 
In the data-driven control setting, however, the order~$n$ of the 
underlying system is unknown, and parametric estimates 
$\widehat{M}_1, \ldots, \widehat{M}_L$ and 
$\widehat{N}_1, \ldots, \widehat{N}_L$, with $\widehat{M}_1 = I_p$ 
are obtained for some choice $L \in \lbrace 1, \ldots, T \rbrace$.
This corresponds to the system identification problem, where
the identified model is given by $\widehat{\mathcal{G}}_{\rm id} = - \widehat{\mathcal{M}}_L^{-1} \widehat{\mathcal{N}}_L$,
and the indirect data-driven control formulation is given by:
\begin{align*}
	\min_{\mathbf{u} \in \real^{mT}} F(\widehat{\mathcal{P}}_{\rm id} \mathbf{u}), 
			\qquad \widehat{\mathcal{P}}_{\rm id} = \left[ \begin{matrix} I \\ - \widehat{\mathcal{M}}_L^{-1} \widehat{\mathcal{N}}_L \end{matrix}	 \right].
\end{align*}
In order to estimate parameters $\lbrace M_k \rbrace_{k=1}^L$ 
and $\lbrace N_k \rbrace_{k=1}^L$ for some choice~$L$,
we reorganize the dataset of input-output behaviors of length~$T$
obtained from control experiments on System~\eqref{eq:LTI_system} into one of
input-output behaviors of length~$L$.
For any length~$T$ behavior $(\mathbf{u}, \mathbf{y})$,
we obtain $T-L+1$ behaviors of length~$L$ as the columns
of $\left[ \begin{matrix} H_L(\mathbf{u})^\top & H_L(\mathbf{y})^\top \end{matrix}	 \right]^\top$,
where for any $\mathbf{z} = (z_1, \ldots, z_T)$, 
$H_L(\mathbf{z})$ is the Hankel matrix of depth~$L$.
For $L \geq n+1$, it follows from~\eqref{eq:IO_delay_rep_vec_form} that:
\begin{align*}
\small
	\left[ \begin{matrix} N_L & \ldots & N_1 & M_{L} & \ldots & M_1 \end{matrix} \right]
				\left[ \begin{matrix} H_L(\mathbf{u}) \\ H_L(\mathbf{y}) - H_L(\mathbf{v})  \end{matrix} \right]
	= 0.
\end{align*}
We henceforth treat the SISO (Single-Input-Single-Output) 
case ($p=m=1$) for the sake of simplicity. We note that our results 
can be readily extended to the MIMO (Multiple-Input-Multiple-Output) case.
For $p=m=1$, the parameters $N_k, M_k$ above are scalars,
and we let $\mathbf{n} = (n_L, \ldots, n_1)$ and $\mathbf{m} = (m_L, \ldots, m_1)$.
Now, we note that the length-$L$ behaviors do not all correspond to
distinct control experiments starting from a zero initial state.
Indeed, any length-$L$ behavior of System~\eqref{eq:LTI_system} satisfies:
\begin{align} \label{eq:length_L_behavior}
\small
\begin{aligned}
	\left( \begin{matrix} \mathbf{u}_{t:t+L-1} \\ \mathbf{y}_{t+1:t+L} \end{matrix} \right) 
					= \left[ \begin{matrix} 0 & I \\ \mathcal{O}_L & \mathcal{G}_L  \end{matrix} \right] 
												\left( \begin{matrix} A x_t \\ \mathbf{u}_{t:t+L-1}	\end{matrix} \right) + \left( \begin{matrix} 0 \\ \mathbf{v}_{t+1:t+L} \end{matrix} \right),
\end{aligned}
\end{align}
where $x_t = \left[ \begin{matrix} A^{t-1}B & \ldots & B \end{matrix} \right] \mathbf{u}_{0:t-1}$,
and the matrices $\mathcal{O}_L$ and $\mathcal{G}_L$ are given by:
\begin{align*}
\footnotesize
	\mathcal{O}_L = \left[ \begin{matrix} C \\  CA \\ \vdots \\  CA^{L-1}  \end{matrix}	 \right], ~
	\mathcal{G}_L =  \left[ \begin{matrix} 	CB & 0 & \ldots & 0 \\
																CAB & CB & \ldots & 0 \\
																	\vdots & \vdots & \ddots & \vdots \\
																		C A^{L-1} B & C A^{L-2} B & \ldots & CB \end{matrix} \right],
\end{align*} 
Assuming that $(A,AB)$ is also a controllable pair, for $L \geq n+1$ we have:
\begin{align} \label{eq:model_based_noise-free}
	\left( \begin{matrix} \mathbf{n}^\top & \mathbf{m}^\top \end{matrix} \right) \left[ \begin{matrix} 0 & I \\ \mathcal{O}_L  & \mathcal{G}_L  \end{matrix} \right]  = 0.
\end{align}
It follows that $\mathbf{n} + \mathcal{G}_L^\top \mathbf{m} = 0$
and $\mathcal{O}_L^\top \mathbf{m} = 0$.
For $k \in \{ 1, \ldots, L-1 \}$, we have $n_k = - \mathcal{G}_L^{(L-k+1) \top} \mathbf{m}$.
Furthermore, since $\mathcal{O}_L^\top \mathbf{m} = 0$,
it follows for $k=L$ that $n_L + \mathcal{G}_L^{(1) \top} \mathbf{m}
= n_L + B^\top \mathcal{O}_L^\top \mathbf{m} =  n_L = 0$.
For the problem of identifying the parameters $(\mathbf{n}, \mathbf{m})$,
we note that we do not have access to the matrices $\mathcal{O}_L$
and $\mathcal{G}_L$ above, nor to the (process) noise $\mathbf{v}$
from experiments. We thereby obtain a data-driven, 
ordinary least squares-based formulation to 
solve~\eqref{eq:model_based_noise-free}:
\begin{align}\label{eq:sys_id_minimization}
\small
\begin{aligned}
 \min_{\tilde{\mathbf{n}}, \tilde{\mathbf{m}} \in \real^L}
																		~&\left( \begin{matrix} \tilde{\mathbf{n}}^\top & \tilde{\mathbf{m}}^\top \end{matrix} \right)
			\frac{1}{\tilde{N}}\sum_{k=1}^N \left[ \begin{matrix} H_L(\mathbf{u}^{(k)}) \\ H_L(\mathbf{y}^{(k)}) \end{matrix} \right] \left[ \begin{matrix} H_L(\mathbf{u}^{(k)}) \\ H_L(\mathbf{y}^{(k)}) \end{matrix} \right]^\top 
																				\left( \begin{matrix} \tilde{\mathbf{n}} \\ \tilde{\mathbf{m}} \end{matrix} \right), \\
																				& \text{s.t.}~~	 \tilde{m}_1 = 1,~ \tilde{n}_L = 0,
\end{aligned}
\end{align}
where $\tilde{N} = N(T-L+1)$.
We now investigate the dependence of the identified parameter mismatch 
on the number of control experiments~$N$, the time horizon~$T$
and the model dimension~$L$.
We first note that:
\begin{align*}
\small
\begin{aligned}
	&\mathbb{E} \left[ \frac{1}{\tilde{N}}
			\sum_{k=1}^N \left[ \begin{matrix} H_L(\mathbf{u}^{(k)}) \\ H_L(\mathbf{y}^{(k)}) \end{matrix} \right] \left[ \begin{matrix} H_L(\mathbf{u}^{(k)}) \\ H_L(\mathbf{y}^{(k)}) \end{matrix} \right]^\top  \right]  \\
	&\qquad = \left[ \begin{matrix}  \Sigma_{\rm u} &  \Sigma_{\rm u} \mathcal{G}_L^\top 
					\\ \mathcal{G}_L \Sigma_{\rm u} & \mathcal{G}_L \Sigma_{\rm u} \mathcal{G}_L^\top  +   \mathcal{O}_L \Sigma_{0}  \mathcal{O}_L^\top + \Sigma_{\rm v} \end{matrix} \right].
\end{aligned}
\end{align*}
Let $\widebar{\mathbf{n}} = \mathbb{E} \left[ \widehat{\mathbf{n}} \right], 
\widebar{\mathbf{m}} = \mathbb{E} \left[ \widehat{\mathbf{m}} \right]$.
It can be shown that since $(\widebar{\mathbf{n}}, \widebar{\mathbf{m}})
= \mathbb{E} \left[ (\widehat{\mathbf{n}}, \widehat{\mathbf{m}}) \right]$,
the pair $(\widebar{\mathbf{n}}, \widebar{\mathbf{m}})$
is the minimizer in~\eqref{eq:sys_id_minimization} 
in expectation. 
For the unconstrained version of~\eqref{eq:sys_id_minimization}, 
in expectation we get:
\begin{align*}
	&\Sigma_{\rm u} \widebar{\mathbf{n}} + \Sigma_{\rm u} \mathcal{G}_L^\top \widebar{\mathbf{m}} = \lambda_{\rm min} \widebar{\mathbf{n}}, \\
	&\mathcal{G}_L \Sigma_{\rm u} \widebar{\mathbf{n}} + \mathcal{G}_L \Sigma_{\rm u} \mathcal{G}_L^\top \widebar{\mathbf{m}} + \mathcal{O}_L \Sigma_{0} \mathcal{O}_L^\top \widebar{\mathbf{m}}
					+	\Sigma_{\rm v} \widebar{\mathbf{m}} = \lambda_{\rm min} \widebar{\mathbf{m}},
\end{align*}
where $\lambda_{\min} \geq 0$ is the smallest eigenvalue of
the covariance matrix above.
It follows that:
\begin{align*}
	\lambda_{\rm min} \mathcal{G}_L \widebar{\mathbf{n}} + \mathcal{O}_L \Sigma_{0} \mathcal{O}_L^\top \widebar{\mathbf{m}}
					+	\Sigma_{\rm v} \widebar{\mathbf{m}} = \lambda_{\rm min} \widebar{\mathbf{m}},
\end{align*}
Note that for $\lambda_{\min} = 0$, we 
will need that $\mathcal{O}_L \Sigma_{0} \mathcal{O}_L^\top + \Sigma_{\rm v}$
has a non-trivial kernel, which is precluded by the fact that
$\Sigma_{\rm v}$ is positive definite and for $L \leq n$,
even $\mathcal{O}_L^\top$ may have a non-trivial kernel.
Thereby, for $\lambda_{\min} >  0$, we have:
\begin{align*}
\footnotesize
\begin{aligned}
	\left( \begin{matrix} \widebar{\mathbf{n}}^\top & \widebar{\mathbf{m}}^\top \end{matrix}  \right)
		\left[ \begin{matrix}	\lambda_{\rm min} \mathcal{G}_L^\top & \Sigma_{\rm u} - \lambda_{\min} I	 \\
													  \left( \mathcal{O}_L \Sigma_{0} \mathcal{O}_L^\top + \Sigma_{\rm v} \right) - \lambda_{\rm min} I 
																			& \mathcal{G}_L \Sigma_{\rm u}		\end{matrix} \right] = 0.
\end{aligned}
\end{align*}
Comparing the above to~\eqref{eq:sys_id_minimization},
we see that in the presence of process noise and/or for 
$L \leq n$, we incur an asymptotic bias
$\left\| (\widebar{\mathbf{n}}, \widebar{\mathbf{m}}) - (\mathbf{n}, \mathbf{m}) \right\| > 0$,
from the fact that $\lambda_{\min} > 0$.

For the non-asymptotic case, we first note that 
the data matrices satisfy:
\begin{align*}
	\left[ \begin{matrix} H_L(\mathbf{u}^{(k)}) \\ H_L(\mathbf{y}^{(k)}) \end{matrix} \right]
	= \left[ \begin{matrix} I \\ \mathcal{G}_L \end{matrix} \right] H_L(\mathbf{u}^{(k)}) 
					+ \left[ \begin{matrix} 0 \\ \mathcal{O}_L \end{matrix} \right] A X_0
							+ \left[ \begin{matrix} 0 \\ H_L(\mathbf{v}^{(k)}) \end{matrix} \right],
\end{align*}
where $X_0$ is the matrix containing the hidden initial states for the
length-$L$ behaviors in the data Hankel matrices.
We see that the sample covariance matrix in~\eqref{eq:sys_id_minimization}
is given by:
\begin{align*}
\footnotesize
\begin{aligned}
\frac{1}{\tilde{N}}\left[ \begin{matrix}
\sum_{k=1}^N H_L(\mathbf{u}^{(k)}) H_L(\mathbf{u}^{(k)})^\top &
\sum_{k=1}^N H_L(\mathbf{u}^{(k)}) H_L(\mathbf{y}^{(k)})^\top \\
\sum_{k=1}^N H_L(\mathbf{y}^{(k)}) H_L(\mathbf{u}^{(k)})^\top &
\sum_{k=1}^N H_L(\mathbf{y}^{(k)}) H_L(\mathbf{y}^{(k)})^\top.
\end{matrix} \right]
\end{aligned}
\end{align*}
We now see that $\sum_{k=1}^N H_L(\mathbf{u}^{(k)}) H_L(\mathbf{u}^{(k)})^\top/\tilde{N} = \Sigma^N_{{\rm uu},L}$, 
where $\Sigma^N_{{\rm uu},L}$ is the covariance of the length-$L$ inputs.
Furthermore, we get $ H_L(\mathbf{u}) H_L(\mathbf{y})^\top /\tilde{N}
= \Sigma^N_{{\rm uu},L} \mathcal{G}_L^\top +  \Sigma^N_{{\rm u x_0},L} A^\top \mathcal{O}_L^\top 
+ \Sigma^N_{{\rm u v},L}$,
where $\Sigma^N_{{\rm u x_0},L}$ and $\Sigma^N_{{\rm u v},L}$
are the sample (cross) covariances between $\mathbf{u}, x_0$ and $\mathbf{u}, \mathbf{v}$ respectively.
A similar expression involving the sample (cross) covariances among 
$\mathbf{u}, x_0, \mathbf{v}$ can be obtained for the
block $\sum_{k=1}^N H_L(\mathbf{y}^{(k)}) H_L(\mathbf{y}^{(k)})^\top /N$.
We note that the cross covariances vanish in expectation
since $\mathbf{u}, x_0, \mathbf{v}$ are independent.
We further note that the matrices $\Sigma^N_{{\rm uu},L},
\Sigma^N_{{\rm u x_0},L}, \Sigma^N_{{\rm u v},L}, \Sigma^N_{{\rm x_0 v},L}$
are $L$-dimensional square matrices, and do not depend on~$T$.
Following a similar analysis as in Theorem~\ref{thm:direct_data-driven},
we can infer that as the dataset size~$N$ increases, they converge to zero with high probability
at a rate $\mathcal{O}\left( \frac{L^2}{N(T-L+1)} \right)$, which is much 
faster than the rate $\mathcal{O}\left( \frac{T^2}{N} \right)$
we obtained for the direct data-driven~case.

\section{A comparison of direct and indirect \\ data-driven control}
\begin{figure*}[!h]
\centering
        \includegraphics[height=0.27\textwidth]{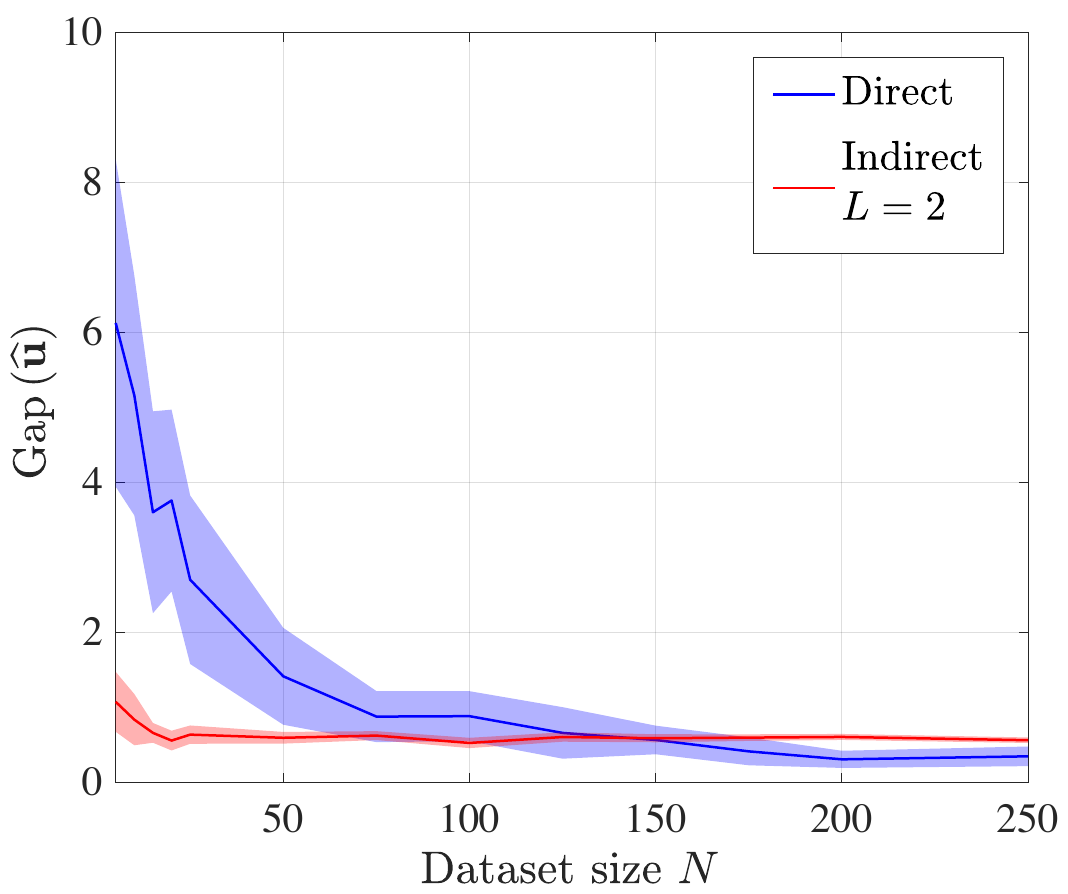}
                \includegraphics[height=0.27\textwidth]{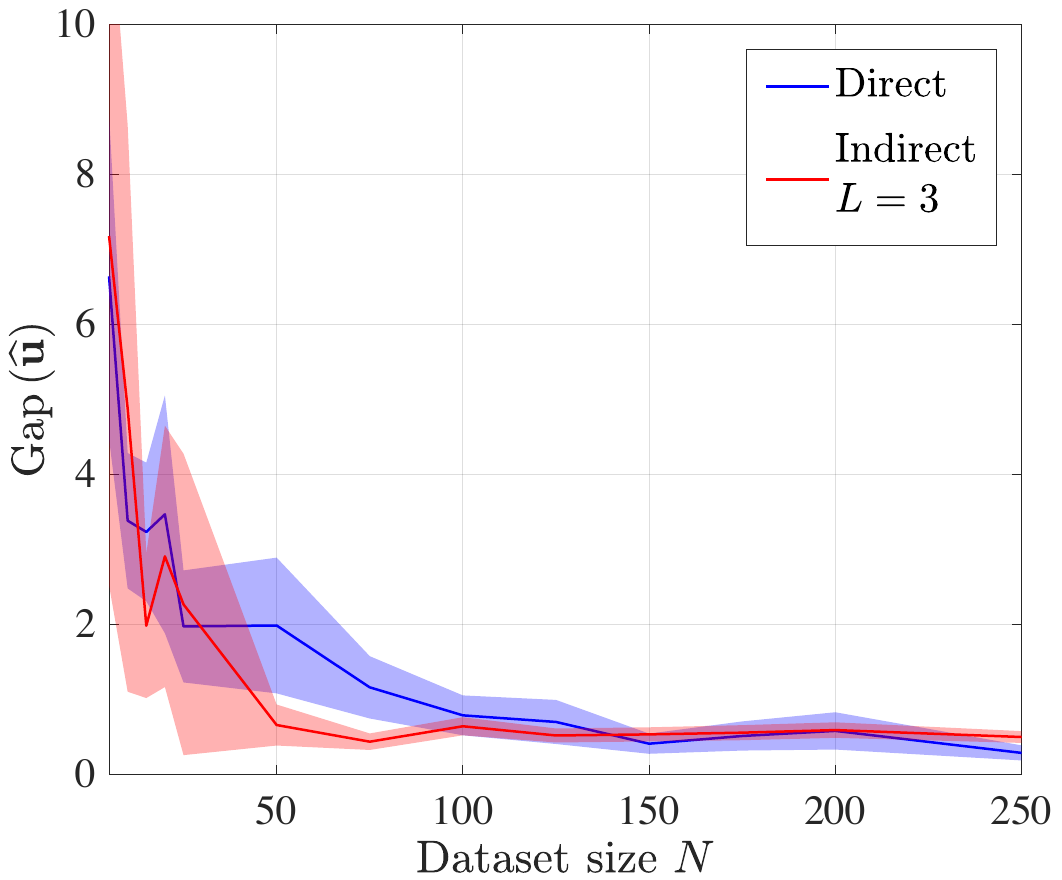}
                        \includegraphics[height=0.27\textwidth]{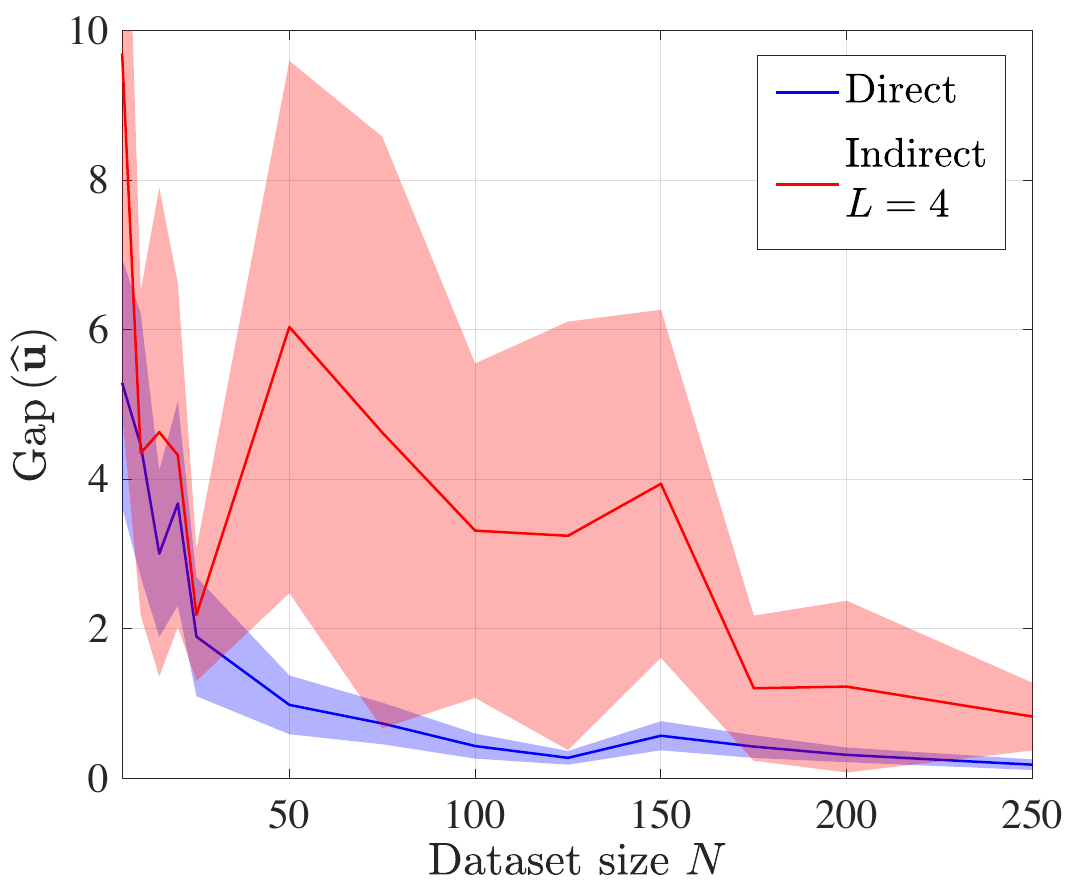}
\caption{The figure shows the dependence of the suboptimality gap 
				${\rm Gap} \left( \widehat{\mathbf{u}} \right)$ defined 
				in~\eqref{eq:suboptimality_gap} on the 
				dataset size~$N$ for direct and indirect data-driven predictive control,
				for three different values $L =2,3,4$ (dimension of the identified system). 
				The underlying system matrices $A,B,C$ in~\eqref{eq:LTI_system} 
				were randomly generated (entries are i.i.d. normal random samples) 
				with $n=3$ and $p=m=1$ (Single Input Single Output case). 
				The task horizon length was chosen to be $T=5$, cost matrices $Q = I$,
				$R = I$ and $y_{\rm ref} = 1$. Inputs in the control experiments are
				i.i.d. samples of the standard normal distribution and the process noise
				are i.i.d. samples of the normal distribution with zero mean and covariance
				$\Omega_{\rm w} = 0.75 I$. The plots show the (empirical) mean suboptimality 
				gap from $50$ trials (sets of control experiments for every $N$), 
				along with the corresponding $95\%$ confidence bands.
				With the direct approach as the reference, we observe that
				the convergence rate of the indirect approach deteriorates as~$L$
				increases, seen as an increase in the empirical mean suboptimality gap
				and a widening of the $95\%$ confidence band, suggesting that the
				indirect approach is prone to overfitting the noise in the dataset 
				at higher values of~$L$.}
		\label{fig:num_expt_fixed_T}
\end{figure*}

\begin{figure*}[!h]
\centering
        \includegraphics[height=0.27\textwidth]{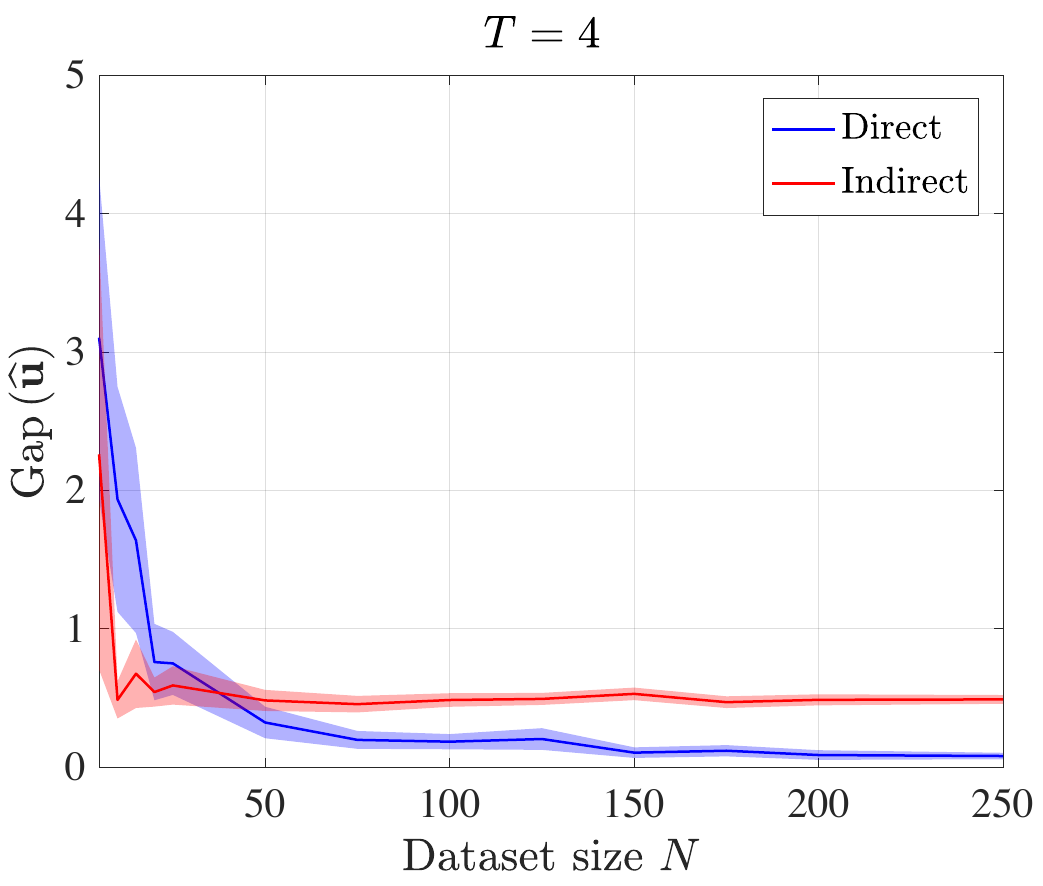}
                \includegraphics[height=0.27\textwidth]{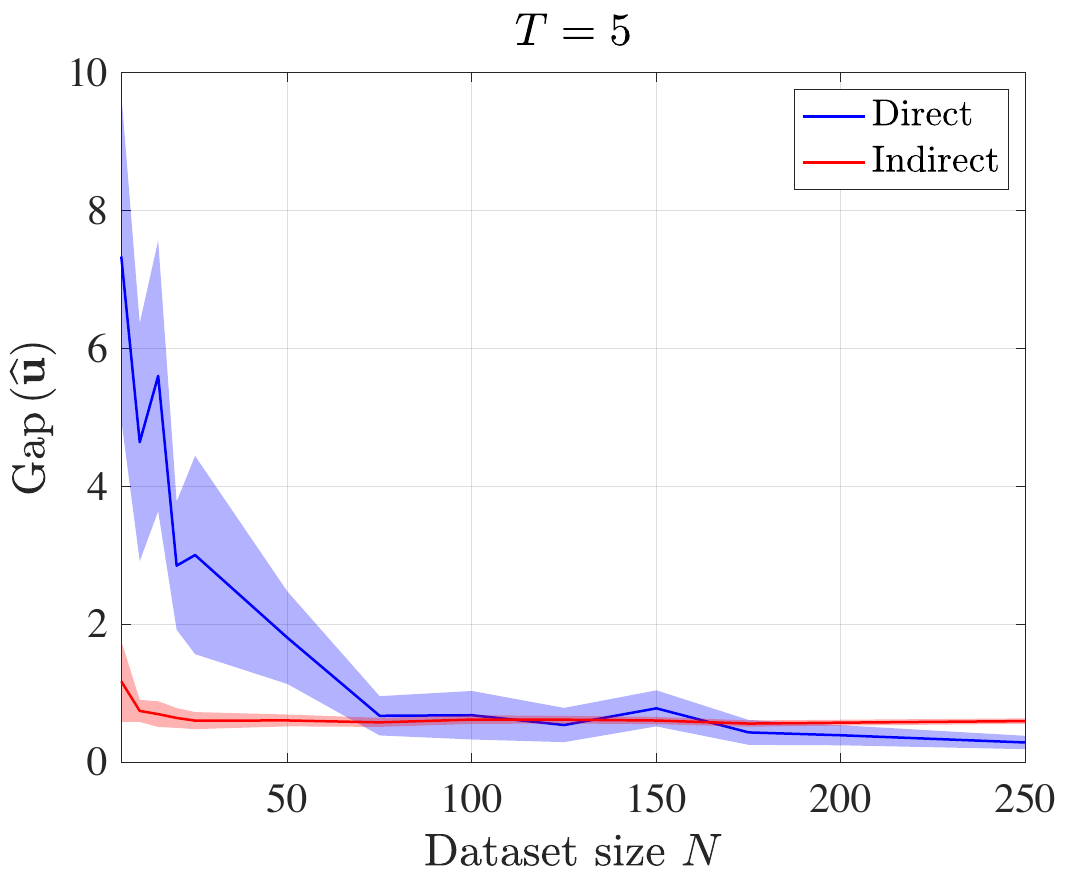}
                        \includegraphics[height=0.27\textwidth]{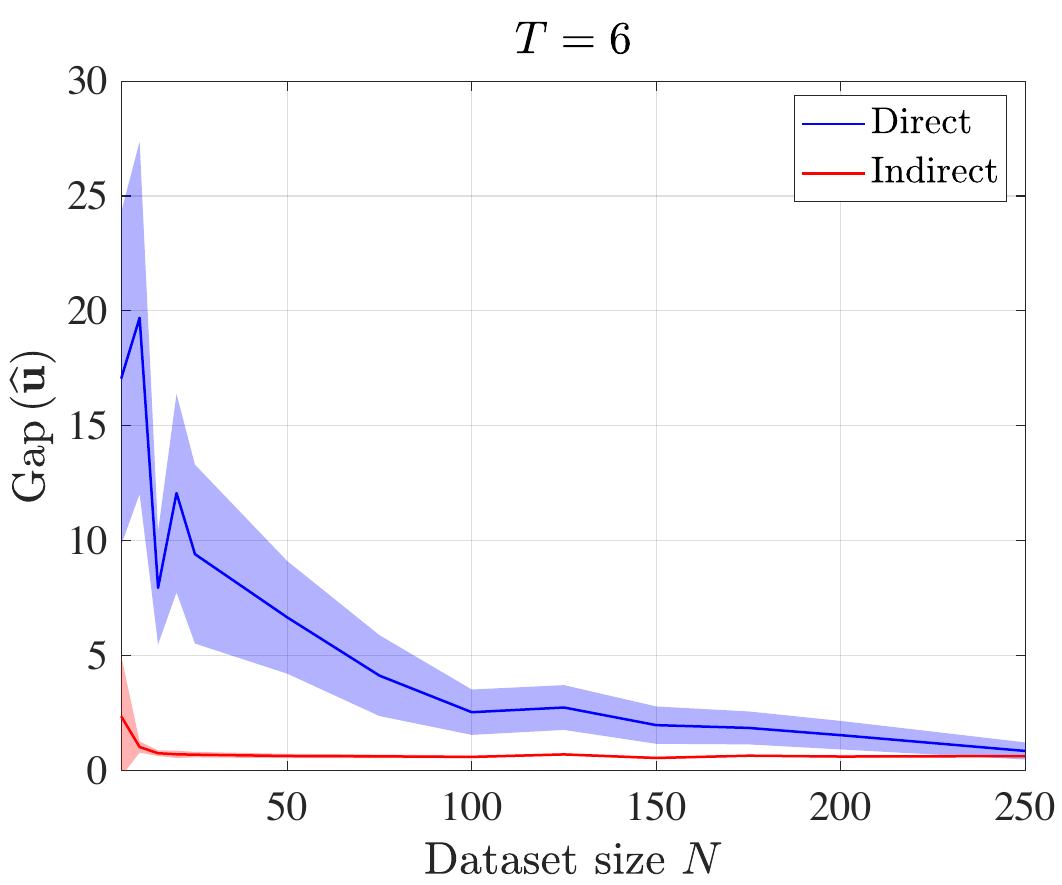}
\caption{The figure shows the dependence of the suboptimality gap 
				${\rm Gap} \left( \widehat{\mathbf{u}} \right)$ defined  
				in~\eqref{eq:suboptimality_gap} on the 
				dataset size~$N$ for direct and indirect data-driven predictive control,
				for three different values $T=4,5,6$ of the control horizon length. The experimental
				setup is otherwise as outlined in Figure~\ref{fig:num_expt_fixed_T},
				with $L=2$ for the indirect~approach. We observe that the sample complexity
				of the direct approach increases with the length $T$ of the control horizon, while
				the indirect approach remains unaffected by the control horizon length.}
		\label{fig:num_expt_variable_T}
\end{figure*}

\begin{figure*}[!h]
\centering
        \includegraphics[height=0.268\textwidth]{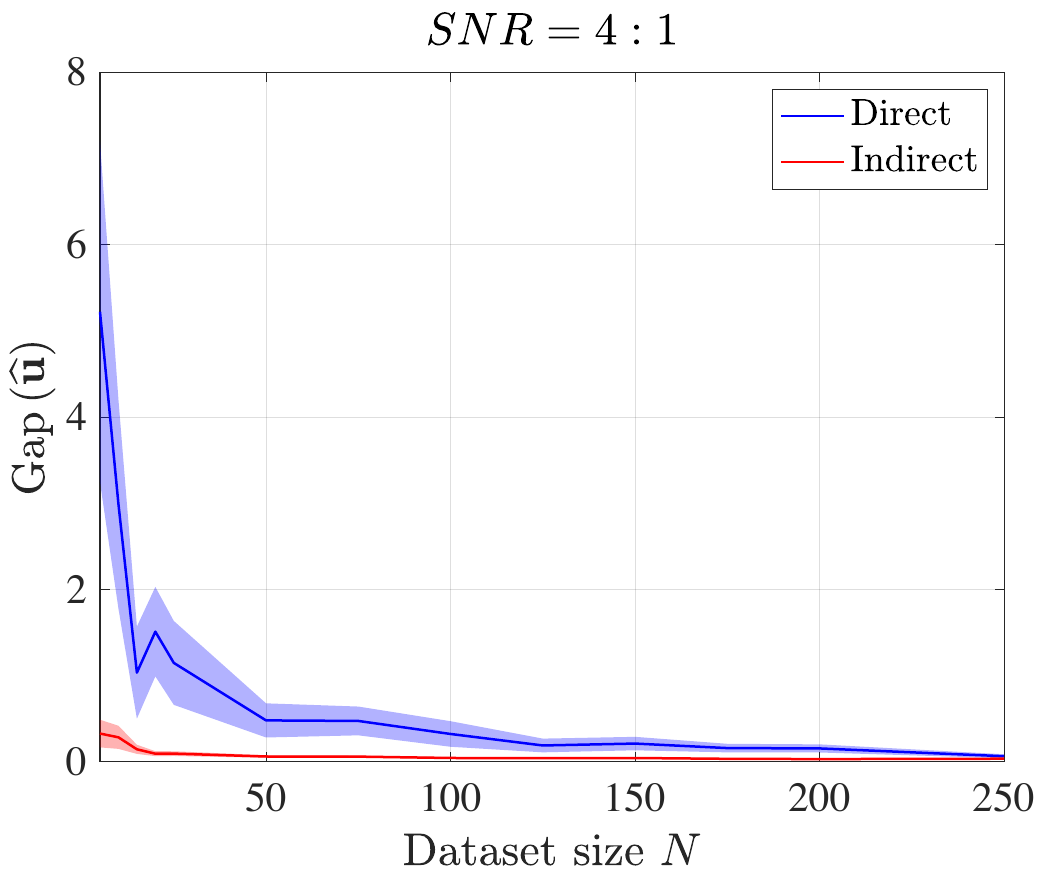}
                \includegraphics[height=0.27\textwidth]{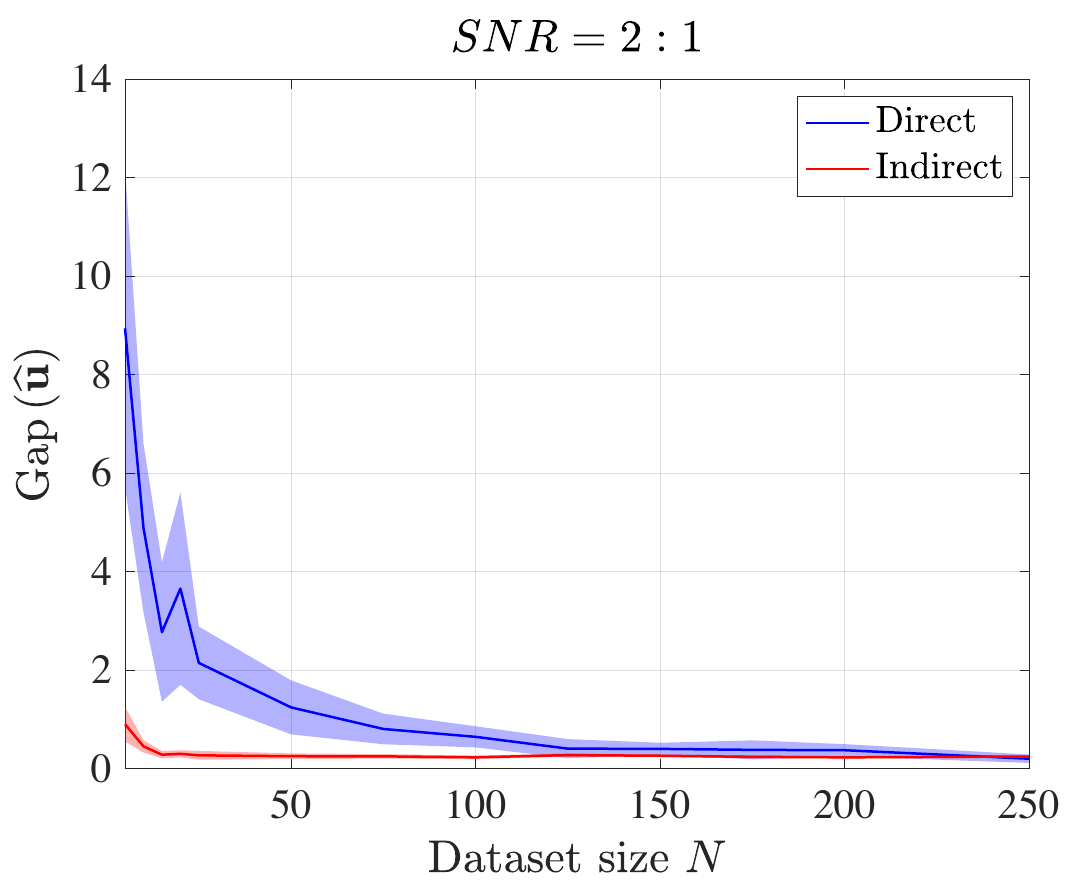}
                        \includegraphics[height=0.27\textwidth]{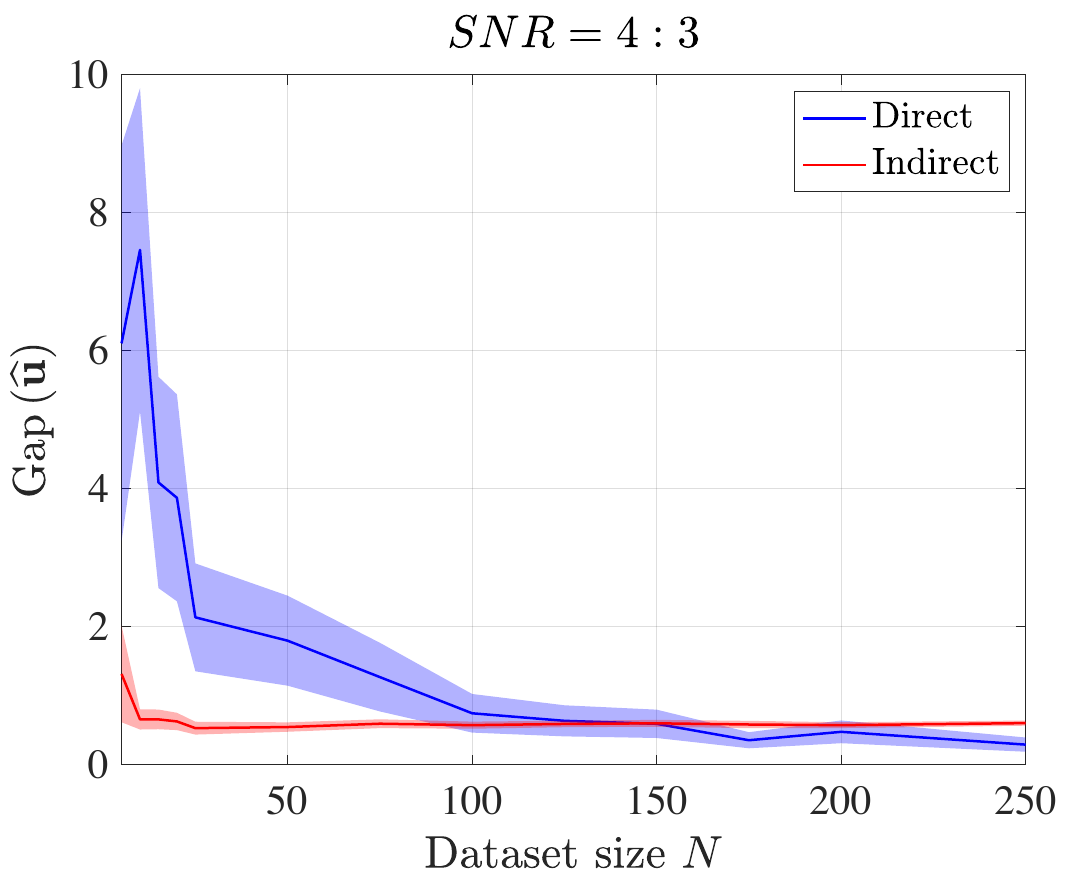}
\caption{The figure shows the dependence of the suboptimality gap 
				${\rm Gap} \left( \widehat{\mathbf{u}} \right)$ defined 
				in~\eqref{eq:suboptimality_gap} on the 
				dataset size~$N$ for direct and indirect data-driven predictive control,
				for three different values of the signal-to-noise ratio (control input to process noise)
				in the control experiments to generate the dataset. 
				The experimental	setup is otherwise as outlined in Figure~\ref{fig:num_expt_fixed_T},
				with length $T=5$ of the control horizon and $L=2$ for the indirect approach.
				We observe that the asymptotic suboptimality gap of the indirect approach
				increases as the signal-to-noise ratio decreases, and subsequently the direct approach outperforms 
				the indirect approach at lower values of~$N$, the dataset size.}
		\label{fig:num_expt_variable_SNR}
\end{figure*}

We now discuss the comparative performance
of the direct and indirect approaches, drawing both from the analytical results
of the previous section and results from numerical experiments (contained
in Figures~\ref{fig:num_expt_fixed_T}, \ref{fig:num_expt_variable_T} and~\ref{fig:num_expt_variable_SNR}).

We see from Theorem~\ref{thm:direct_data-driven} that 
the implicit model error in the direct approach vanishes asymptotically
as $N \rightarrow \infty$, and it thereby follows that the suboptimality gap vanishes 
asymptotically with the size of the dataset. However, we also see that the direct approach
has high sample complexity, resulting in potentially large variance for finite values of~$N$,
as clearly seen from the wider $95\%$ confidence bands for the
direct approach in Figure~\ref{fig:num_expt_variable_SNR}.

From the discussion in the previous section and 
Figure~\ref{fig:num_expt_variable_SNR}, we see 
that the ordinary least squares-based indirect approach 
incurs an asymptotic bias (and consequently, a suboptimality gap) 
due to the process noise covariance and
an incorrect assumption on the dimension
of the underlying system. While the asymptotic bias due to 
process noise can be partially mitigated by a total least
squares-based approach~\cite{IM-SVH:07} for high signal-to-noise ratio
regimes, the bias due to an incorrect assumption on the 
model dimension will remain. Moreover, we note that a 
total least squares-based approach has pitfalls in low 
signal-to-noise ratio regimes in comparison to ordinary
least squares, and will potentially have high variance
(where ordinary least squares has a distinct
advantage over the direct approach), 
and we believe that there exists no invariably
superior candidate for the indirect approach.

We note that the indirect approach has lower sample complexity than the direct
approach owing to the intermediate step of identifying a lower $L$-dimensional
model from the data. Furthermore, we observe from experiments that the indirect approach 
with a large $L$ ($\geq n+1$) overfits 
noise in the data resulting in high variance (Figure~\ref{fig:num_expt_fixed_T}),
and the direct approach outperforms the indirect approach in this case. 
We also see from Figure~\ref{fig:num_expt_fixed_T} that 
the indirect approach incurs a larger asymptotic bias with a lower value of $L$ ($\leq n$) 
but has lower variance for finite values of~$N$.

Since the direct approach relies on
an implicit $T$-dimensional model of the input-output 
behavior, we see that the indirect approach mitigates 
some of the hurdles faced by the (high dimensional) 
direct approach by explicitly identifying a lower 
$L$-dimensional model. 
Intuitively, we see that the implicit $T$-dimensional model 
in the direct approach effectively contains $\mathcal{O}\left( T^2 \right)$
parameters, while the identified $L$-dimensional model
in the indirect approach contains $\mathcal{O}\left( L \right)$
parameters. This explains the poor scaling behavior of the direct approach 
observed in Figure~\ref{fig:num_expt_variable_T}, where the performance
of the direct approach deteriorates with increasing~$T$,
while the indirect approach remains unaffected (since $L$
is not a function of~$T$).
Furthermore, identifying a lower $L$-dimensional 
model not only has lower sample complexity, but also 
results in the expansion of the dataset available for
identification. This can be seen from the fact that
the dataset of $T$-dimensional behaviors
contains~$N$ samples, whereas by constructing 
$L$-dimensional Hankel matrices from $T$-dimensional
behaviors essentially results in a dataset containing
$N(T-L+1)$ samples (it must however be noted that
the $L$-dimensional behavior samples are not all~i.i.d.).

\section{Concluding remarks}
Our study reveals the existence of two non-asymptotic regimes 
for the performance of direct and indirect data-driven predictive
control designs, which precludes any conclusion that either
of the two approaches is invariably superior. It also
suggests that interpolating between the direct and indirect approaches 
may be useful in managing the underlying tradeoffs effectively in practice,
a direction that has recently been pursued in~\cite{FD-JC-IM:21}.
We note that for our comparison we chose a
certainty-equivalent control design with
ordinary least squares-based identification 
as the candidate for the indirect approach. 
While our results are, in principle, restricted 
by this particular choice,
we believe that the qualitative insights we obtain 
on the comparative performances
are much more general and not constrained by
this choice. However, we do not comment on the 
capabilities of potentially more sophisticated 
techniques (for both direct and indirect data-driven control) 
in better managing the underlying tradeoffs,
and a more fundamental analysis is likely necessary for
a truly technique-agnostic comparison.
Furthermore, a fine-grained analysis of the implications for
data-driven control of the phenomena 
of measure concentration, data sparsity and 
scaling of signal-to-noise ratio may afford a finer 
characterization of the non-asymptotic regimes of performance
and result in practically useful heuristics.

\bibliographystyle{unsrt}
\bibliography{alias,Main,FP,New}

\appendix
\section{}
\subsection{Proof of Lemma~\ref{lemma:certainty_equiv_model-based}} \label{app:proof_certainty_equiv}
By convexity of $F$ and Jensen's inequality, we have:
\begin{align*}
	\mathbb{E}_{\mathbf{v}} \left[ F(\mathbf{u}, \mathcal{G} \mathbf{u} + \mathbf{v}) \right] &\geq F( \mathbf{u} , \mathbb{E}_{\mathbf{v}} \left[ \mathcal{G} \mathbf{u} + \mathbf{v} \right]) \\
																	  &= F(\mathbf{u}, \mathcal{G} \mathbf{u}) \\
																	  &\geq \min_{\mathbf{u}} F(\mathbf{u}, \mathcal{G} \mathbf{u}).
\end{align*}
Let $G(\mathbf{u}) = \mathbb{E}_{\mathbf{v}} \left[ F(\mathbf{u}, \mathcal{G} \mathbf{u} + \mathbf{v}) \right]$. 
For any $\mathbf{u}_1$ and $\mathbf{u}_2$, we have:
\begin{align*}
	G(\mathbf{u}_2) - G(\mathbf{u}_1) 
						&= \mathbb{E}_{\mathbf{v}} \left[ F(\mathbf{u}_2, \mathcal{G} \mathbf{u}_2 + \mathbf{v}) - F(\mathbf{u}_1, \mathcal{G} \mathbf{u}_1 + \mathbf{v}) \right] \\
						&\geq \mathbb{E}_{\mathbf{v}} \left[ \nabla^\top F (\mathbf{u}_1, \mathcal{G} \mathbf{u}_1 + \mathbf{v}) 
																	\left( \begin{matrix} \mathbf{u}_2 - \mathbf{u}_1 \\ \mathcal{G} (\mathbf{u}_2 - \mathbf{u}_1) \end{matrix} \right) \right] \\
						&= \mathbb{E}_{\mathbf{v}} \left[ \nabla^\top F (\mathbf{u}_1, \mathcal{G} \mathbf{u}_1 + \mathbf{v}) 
													\left[ \begin{matrix} I \\ \mathcal{G} \end{matrix} \right]  \right]  (\mathbf{u}_2 - \mathbf{u}_1) \\
						&= \mathbb{E}_{\mathbf{v}} \left[ \nabla_{\mathbf{u}}^\top F (\mathbf{u}_1, \mathcal{G} \mathbf{u}_1 + \mathbf{v}) \right] (\mathbf{u}_2 - \mathbf{u}_1).
\end{align*}
By the Dominated Convergence Theorem, we get:
\begin{align*}
	\mathbb{E}_{\mathbf{v}} \left[ \nabla_{\mathbf{u}} F (\mathbf{u}, \mathcal{G} \mathbf{u} + \mathbf{v}) \right]
			= \nabla_{\mathbf{u}} \mathbb{E}_{\mathbf{v}} \left[ F (\mathbf{u}, \mathcal{G} \mathbf{u} + \mathbf{v}) \right]
				= \nabla_{\mathbf{u}} G(\mathbf{u}).
\end{align*}
Substituting in the inequality above, we get:
\begin{align*}
	G(\mathbf{u}_2) - G(\mathbf{u}_1)  \geq \nabla_{\mathbf{u}}^\top G(\mathbf{u}) (\mathbf{u}_2 - \mathbf{u}_1).
\end{align*}
This establishes that the function $G$ is also convex.
Furthermore, since $F$ is quadratic, it follows that $\nabla F$ is affine.
Now, for $\mathbf{u}^* = \arg \min_{\mathbf{u}}  F(\mathbf{u}, \mathcal{G} \mathbf{u})$, we have
$\nabla_{\mathbf{u}} G(\mathbf{u}^*) 
= \mathbb{E}_{\mathbf{v}} \left[ \nabla_{\mathbf{u}} F (\mathbf{u}^*, \mathcal{G} \mathbf{u}^* + \mathbf{v}) \right]
=  \nabla_{\mathbf{u}} F (\mathbf{u}^*, \mathcal{G} \mathbf{u}^*) 
= 0$ (by strict convexity of $F$ it also follows that $\nabla_{\mathbf{u}} G(\mathbf{u}) = 0$
only if $\mathbf{u} = \mathbf{u}^*$). Therefore, for any $\mathbf{u} \in \real^{mT}$, we get that 
$G(\mathbf{u}) - G(\mathbf{u}^*) \geq \nabla_{\mathbf{u}}^\top G(\mathbf{u}^*) (\mathbf{u} - \mathbf{u}^*)
 = 0$, which implies that $\mathbf{u}^*$ is the global minimizer
 of $G$. Thus, we have 
$\mathbf{u}^* = \arg \min_{\mathbf{u}} \mathbb{E}_{\mathbf{v}} \left[ F(\mathbf{u}, \mathcal{G} \mathbf{u} + \mathbf{v}) \right]$.

\subsection{Proof of Theorem~\ref{thm:direct_data-driven}} \label{app:proof_direct}
\noindent The implicit model error $\Delta_{\rm direct}$ is given by:
\begin{align*}
	&\Delta_{\rm direct} = V U^\dagger = V U^\top \left( U U^\top \right)^{-1} \\
	&= \sum_{k=1}^N \mathbf{v}^{(k)} \mathbf{u}^{(k) \top} \left( \sum_{k=1}^N \mathbf{u}^{(k)} \mathbf{u}^{(k) \top} \right)^{-1} \\
	&=  \sum_{k=1}^N \mathcal{G}' \mathbf{w}^{(k)} \mathbf{u}^{(k)\top} \left( \sum_{k=1}^N \mathbf{u}^{(k)} \mathbf{u}^{(k) \top} \right)^{-1}.
\end{align*}
Since $\mathbf{w}$ and $\mathbf{u}$ are independent and 
$\mathbb{E}[\mathbf{w}] = 0$, we have:
\begin{align*}
\footnotesize
\begin{aligned}
	&\mathbb{E} \left[ \Delta_{\rm direct} \right] = \sum_{k=1}^N \mathcal{G}' \mathbb{E} \left[ \mathbf{w}^{(k)} \right] \cdot
															\mathbb{E} \left[ \mathbf{u}^{(k)^\top} \left( \sum_{k=1}^N \mathbf{u}^{(k)} \mathbf{u}^{(k) \top} \right)^{-1} \right] \\
								&= 0.
\end{aligned}
\end{align*}
Let $\Sigma_{\rm uu}^N = \sum_{k=1}^N \mathbf{u}^{(k)} \mathbf{u}^{(k) \top} / N$
be the sample covariance of $\mathbf{u}$. 
Following the reasoning above, we also get that:
\begin{align*}
\small
\begin{aligned}
	\mathbb{E} \left[ \frac{1}{N} V U^\top \right]
		= \mathbb{E} \left[ \frac{1}{N} \sum_{k=1}^N \mathcal{G}' \mathbf{w}^{(k)} \mathbf{u}^{(k) \top} \right] = 0.
\end{aligned}
\end{align*}
Now, we have:
\begin{align*}
\footnotesize
\begin{aligned}
	\left\| \Delta_{\rm direct} \right\|_F \leq \left\| \frac{1}{N} V U^\top \right\|_F \left\| {\Sigma_{\rm uu}^N}^{-1} \right\|_2 
																= \frac{1}{\sigma_{\min} \left( \Sigma_{\rm uu}^N \right)} \left\| \frac{1}{N} V U^\top \right\|_F.
\end{aligned}
\end{align*}
Furthermore, we have: 
\begin{align*}
\small
\begin{aligned}
	\left\| \frac{1}{N} V U^\top \right\|_F  
				= \sqrt{\sum_{i,j=1}^T \left| \left( \frac{1}{N}  VU^\top \right)_{ij} \right|^2 }
				\leq T \left\| \frac{1}{N} VU^\top \right\|_{\max},
\end{aligned}
\end{align*}
and it follows that:
\begin{align*}
	\left\| \Delta_{\rm direct} \right\|_F \leq \frac{T}{\sigma_{\min} \left( \Sigma_{\rm uu}^N \right)} \left\| \frac{1}{N} VU^\top \right\|_{\max}.
\end{align*}
We now have:
\begin{align*}
	&\mathbb{P} \left\{ \left\| \Delta_{\rm direct} \right\|_F \geq \epsilon \right\}
				\leq \mathbb{P} \left\{ \frac{T}{\sigma_{\min} \left( \Sigma_{\rm uu}^N \right)} \left\| \frac{1}{N} VU^\top \right\|_{\max} \geq \epsilon \right\} \\
				&= \mathbb{P} \left\{ \left\| \frac{1}{N} VU^\top \right\|_{\max} \geq  \frac{\sigma_{\min} \left( \Sigma_{\rm uu}^N \right)}{T}  \epsilon \right\}.
\end{align*}
Applying Chebyshev's inequality, we obtain:
\begin{align*}
\small
\begin{aligned}
	&\mathbb{P} \left \lbrace \left| \frac{1}{N} \left( VU^\top \right)_{ij} \right| \geq \delta \right \rbrace \\
	&= \mathbb{P} \left \lbrace \left| \frac{1}{N} \left( VU^\top \right)_{ij} - \mathbb{E} \left[ \frac{1}{N} \left( VU^\top \right)_{ij} \right] \right| \geq \delta \right \rbrace \\
	&\leq \frac{1}{\delta^2} \mathrm{Var} \left( \frac{1}{N} \left( VU^\top \right)_{ij} \right)
		= \frac{1}{\delta^2} \mathrm{Var} \left( \frac{1}{N} \sum_{k=1}^N 	v_i^{(k)} u_j^{(k)} \right) \\
	&= \frac{1}{\delta^2 N^2} \sum_{k=1}^N \mathrm{Var} \left( v_i^{(k)} u_j^{(k)} \right) \\
	&= \frac{1}{\delta^2 N^2} \sum_{k=1}^N \mathrm{Var} \left( v_i^{(k)} \right)  \mathrm{Var} \left( u_j^{(k)} \right) \\
	&= \frac{1}{\delta^2 N}	\mathrm{Var} \left( v_i \right)  \mathrm{Var} \left( u_j \right).
\end{aligned}
\end{align*}
We note that for all $i \in \{ 1, \ldots, T \}$, 
$\mathrm{Var} \left( v_i \right) \leq \mathrm{Var} \left( v_T \right)
= \mathrm{Var} \left( \sum_{t=0}^{T-1} CA^t w_t \right)
= \sum_{t=0}^{T-1} \mathrm{Var} \left( CA^t w_t \right)
= \sum_{t=0}^{T-1}  CA^t \Omega_{\rm w} A^{t \top} C^\top
= \sigma_{\rm w}$,
and $\mathrm{Var} \left( u_j \right) = \sigma_{\rm u}$,
and we get $\mathbb{P} \left \{ \left\| VU^\top / N \right\|_{\max} \geq \delta \right\} 
\leq \sigma_{\rm w}  \sigma_{\rm u} / \delta^2 N$.
and it follows that:
\begin{align*}
\begin{aligned}
	\mathbb{P} \left\{ \left\| \Delta_{\rm direct} \right\|_F \geq \epsilon \right\}
		\leq \frac{T^2}{N \epsilon^2} \cdot \frac{\sigma_{\rm w} \sigma_{\rm u}}{\sigma_{\min}^2 \left( \Sigma_{\rm uu}^N \right)}.
\end{aligned}
\end{align*}

\end{document}